# Joint Optimization of Transfer Location and Capacity in a Multimodal Transport Network: Bilevel Modeling and Paradoxes

Yu Jiang [a, *], Ye Jiao [b,c], Chen Jun [b,c]

[a] DTU Transport, Department of Technology, Management and Economics
Technical University of Denmark, Denmark, 2800
[b] School of Transportation
Southeast University, Nanjing, Jiangsu, 210096
[c] Jiangsu Province Collaborative Innovation Center of Modern Urban Traffic Technologies
Southeast University, Nanjing, Jiangsu, 210096

**ABSTRACT**
With the growing attention towards developing the multimodal transport system to enhance urban mobility, there is an increasing need to construct new, rebuild or expand existing infrastructure to facilitate existing and accommodate newly generated travel demand. Therefore, this paper develops a bilevel model to simultaneously determine the location and capacity of the transfer infrastructure to be built considering elastic demand in a multimodal transport network. The upper level problem is formulated as a mixed integer linear programming problem, while the lower level problem is the capacitated combined trip distribution/assignment model that depicts both destination and route choices of travelers via the multinomial logit formula. To solve the model, the paper develops a matheuristics algorithm that integrates a Genetic Algorithm and a successive linear programming solution approach. Numerical studies are conducted to demonstrate the existence and examine two Braess-like paradox phenomena in a multimodal transport network. The first one states that under fixed demand constructing parking spaces to stimulate the usage of Park-and-Ride service could deteriorate the system performance, measured by the total passengers' travel time, while the second one reveals that under variable demand increasing the parking capacity for the Park-and-Ride services to promote the usages may fail, represented by the decline in its modal share. Meanwhile, the last experiment suggests that constructing transfer infrastructures at distributed stations outperforms building a large transfer center in terms of attracting travelers using sustainable transit modes.

[*] Corresponding author: yujiang@dtu.dk
ORCID: 0000-0001-9461-633X

**INTRODUCTION**

With the trend of urban sprawl and the growth of house pricing, more and more people choose to live far from the city center. Thus, the necessity of providing sufficient mobility for these people has attracted more and more attention from both industry and academia. A well-recognized appealing solution is the provision of integrated multimodal transport services. The methodology involved in providing efficient multimodal services can be classified as the multimodal network design problem.

One of the key components in the multimodal network design problem is to define the decision variables in a multimodal network and how these variables affect traveler' behavior, such as mode and route choices. In general, existing studies can be classified into three approaches. In the first approach, the decision variables associated with each transport mode, i.e., bus route, frequency, or link capacity, are exclusive to each transport mode. Therefore, traffic/transit assignment model can be conducted independently on each transport subnetwork by employing a symmetric link performance function, i.e., BPR function (Lee and Vuchic, 2005; Beltran et al., 2009; Szeto et al., 2010). This approach ignores the fact that the decision variables of one transport mode could af fect the network topology of another transport mode and different transport mode could compete for the same limited infrastructure resources. For example, when buses and cars share the same lane, the construction or expanding a bus lane leads to the capacity reduction for cars. Such an issue is addressed by the second approach which designs the allocation of exclusive lanes to specific transportation modes (Elshafei, 2006; Mesbah et al., 2008; Li and Ju, 2009; Yao et al., 2015). The above two approaches assume that travelers' mode and route choice dedicate to one specific transport mode. In another word, the traveler' trips only contain one transport mode[†]. This overlooks the intermodal travel behavior, meaning that a traveler utilizes more than one mode during a trip, such as Park-and-Ride (P+R) and Bike-and-Ride (B+R). Such intermodal travel is becoming a prevailing option for commuters. In view of this, the third approach emerges, which designs the transfer location in a multimodal transport network (Arnold et al., 2004; Alumur et al., 2012). This approach captures a more realistic travel behavior, where intermodal trips are considered in the route choices. Nevertheless, to our best knowledge, the studies within this approach are limited and there is no existing study on designing the location and capacity of transfer infrastructure simultaneously.

Irrespective of the modelling approaches, a common framework for formulating the multimodal network design problem is the bilevel programming. It has been widely adopted in solving the network design problem for traffic network (Yang et al., 2000; Huang et al., 2001; Szeto et al., 2015; Jiang and Szeto, 2015) or public transport network (Gao et al., 2005; Fan and Machemehl, 2011; Szeto et al., 2011; Yan et al., 2013). In the context of the multimodal network design, a lot of efforts have been paid for during the last decade. Fan et al. (2013) developed a bilevel programming model for locating P+R facilities considering the interaction between decision-makers and commuters via incorporating a stochastic user equilibrium model as their lower level problem. Yu et al. (2015) devised a bilevel model for the bus lane distribution problem, where their upper level problem is to minimize the average travel time of travelers considering the balance on the service levels among different modes and the lower level problem is a multimodal transport network equilibrium problem. Recently, Huang et al. (2018) proposed a bilevel programming model to reconfigure the bus services in favor of integrating multimodal transit

---

[†] Walk is considered to access or egress a mode instead of an independent mode to complete a trip

system. The objectives of their upper level and lower level models are respectively to minimize the sum of the total travel time and operation cost and the sum of the total in-vehicle time and total waiting time between all OD pairs.

To solve a bilevel programming problem, most existing studies adopt heuristic (Yang, 1995; Chiou, 2005; Angelo and Barbosa, 2015), metaheuristic (Koh, 2007; Fan et al., 2013), or matheuristics (Szeto and Jiang, 2012; Szeto and Jiang, 2014; Carosi et al., 2019; Liu et al., 2019) to solve it given that a bilevel network design problem is a well-known NP-hard problem. In this study, we develop a matheuristics algorithm to solve our model. The algorithm relies on a Genetic Algorithm to generalize variables associated with transfer locations and capacities and solve the resultant model with a successful linear approximation model developed by Yang et al. (2000).

Other than developing a mathematical model for the multimodal network design problem, this study also reveals a Braess-like paradox phenomenon in the multimodal network problem. Although the Braess paradox has been in-depth examined in both traffic and transit network, to our best knowledge, there is no such observations in a multimodal transport network. Therefore, this study is motived to demonstrate its existence.

To sum up, the contributions of this study include:
1. Developing a bilevel model for the multimodal network design problem to determine the transfer location and capacity simultaneously.
2. Developing a matheuristic solution algorithm that combines the Genetic Algorithm and a successful linear approximation method to solve the model.
3. Demonstrating the existence of a Brass-like paradox phenomenon in a multimodal transport network
4. Conducting various experiments to examine the property of the proposed model as well as the paradox phenomenon

The remaining of the paper is organized as follows: Section 2 introduces the super network representation for the multimodal transport network, assumptions and the notations, and then present the bilevel formulation. Afterwards, the solution method is developed in Section 3. Section 4 depicts the experiments to illustrate the existence of the Brass-like paradox in a multimodal transport network. Finally, Section 5 gives the conclusions and points out future directions.

## FORMULATION

### Network Representation, Problem Description, and Assumptions

We consider an urban transport system and adopt a supernetwork (Sheffi, 1995) approach to depict the multimodal transportation networks. The supernetwork is denoted by $G = (N, A)$, where $N$ and $A$ represent the set of nodes and links, respectively. It contains $|M|$ subnetworks and each subnetwork, represented by $G^m = (N^m, A^m)$, corresponds to a travel mode $m$, i.e., auto, bus, or bike. A commuter traveling between nodes $o$ and $d$ could travel via either a single mode in one subnetwork, i.e., metro, bus, car, or bicycle, or an intermodal trip through multiple subnetworks, i.e., P+R, B+R, etc. Accordingly, we use sets $M$ and $\bar{M}$ to denote the travel mode that utilize one subnetwork and multiple subnetworks, respectively. For travel mode $m \in \bar{M}$, the set of potential transfer nodes, through which travelers can transfer from one subnetwork to another, is dented as $\bar{N}^m$ and the set of dummy transfer links connecting the transfer nodes between two subnetworks are denoted by $\bar{A}$.

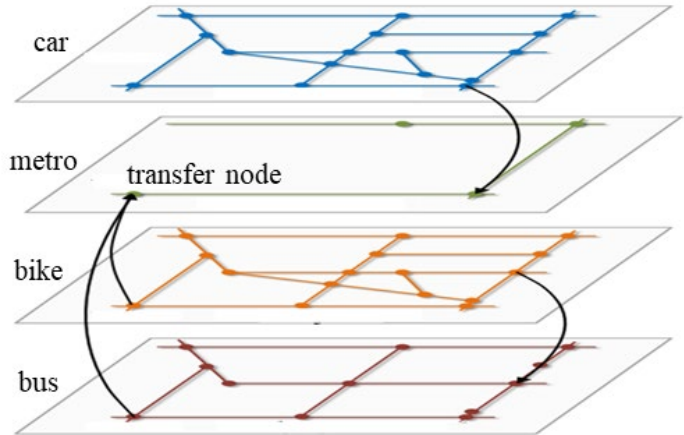

**Figure 1** Multimodal Transportation Network Representation

The multimodal network design problem considered in this paper is twofold. 1) Select locations in the network to construct infrastructures so that an intermodal trip can transfer through the selected locations; 2) For the selected locations, design the infrastructure capacity to accommodate transfer flow. In line with the literature in the multimodal network design problem, the following assumptions are made: A1) The link performance function on each subnetwork is independent. This is reasonable between metro and road networks, and acceptable for car and bus networks by assuming that buses are operated in exclusive bus lanes; A2) For simplicity, the link performance function is continuous and differential and is a function of the flow travelling on the link; A3) The flow on auto subnetwork links refers to vehicle flow, while the flow on metro, bicycle and bus links refer to passenger flow. The vehicular flow are transformed into the passenger flow by the car occupancy rate; A4) The changes in the network affect residents' travel behavior including trip destination, travel mode, and routes (Yang et al., 2000), and it is assumed that the commuters' mode choice and route choice behavior follow the multinomial logit distribution in line with literature; A5) For links on each subnetwork $a \in A^m$, the soft capacity constraint is imposed, meaning that the flow is allowed to be greater than the capacity at additional cost. For the transfer links $a \in \bar{A}$, the hard capacity constraint is imposed, presenting the limited transfer capacity, such as limited parking spaces in the P+R mode or limited docking slots for B+R mode.

**Notations**

The following notations are used in the proposed bilevel model.

*Sets*

| | |
|---|---|
| $M$, $\bar{M}$ | set of travel mode utilizing single and multiple transport modes, respectively |
| $A^m$ | set of links on subnetwork $m \in M$ |
| $\bar{N}^m$ | set of transfer nodes of mode $m \in \bar{M}$ |
| $\bar{A}^m$ | set of transfer links of mode $m \in \bar{M}$ |
| $R, S$ | set of origins and destinations, respectively |
| $P_{rs}^m$ | set of paths connecting OD pair $rs$ via mode $m$ |

*Indices*

| | |
|---|---|
| $m$ | index of travel mode |
| $n$ | index of node |
| $p, r, s$ | index of path, origin, and destination |

*Parameters*

| | |
|---|---|
| $o_r^{\max}, d_s^{\max}$ | maximum demand can be generated at origin $r$ and destination $s$, respectively |
| $o_r^0, d_s^0$ | existing travel demand generated at origin $r$ and destination $d$, respectively |
| $q_{rs}^0$ | existing travel demand between nodes $r$ and $s$ |
| $C_a$ | capacity associated with link $a$ |
| $B$ | total budget |
| $C_n^{m,\max}$ | maximum transfer capacity could be built at node $n$ for mode $m \in \bar{M}$ |
| $C_n^{m,\min}$ | minimum transfer capacity should be built at node $n$ for mode $m \in \bar{M}$ |

*Variables and function of variables*

| | |
|---|---|
| $\bar{c}_a^m$ | designed capacity for mode $m \in \bar{M}$ at a transfer link $a \in \bar{A}_m$ |
| $\bar{v}_n^m$ | transfer flow at node $n$ of mode $m \in \bar{M}$ |
| $\bar{t}_n^m$ | transfer travel time at node $n$ of mode $m \in \bar{M}$ |
| $\delta_{ap}$ | $=1$, if link $a$ is on path $p$, otherwise it equals 0 |
| $\Delta_{np}$ | $=1$, if node $n$ is on path $p$, otherwise it equals 0 |
| $\xi_n^m$ | $=1$, if node $n$ is selected as a transfer location for mode $m \in \bar{M}$, otherwise it equals 0 |
| $v_a^m$ | flow travelling on link $a$ of mode $m \in M$ |
| $f_p^{m0}$ | flow travelling via path $p$ of mode $m$ under existing demand |
| $f_p^{m+}$ | flow travelling via path $p$ of mode $m$ obtained from newly generated demand |
| $o_r^+, d_s^+$ | newly generated travel demand at origin $r$ and destination $d$, respectively |
| $q_{rs}^+$ | newly generated travel demand between nodes $r$ and $s$ |
| $G(\xi)$ | cost function for constructing transfer infrastructure |

**Bilevel Formulation**

A bilevel optimization model is developed for the multimodal network design problem and presented in the following two subsections.

*Upper Level Problem*

$$\max_{\xi,c,q} z = \sum_{r \in R} \sum_{s \in S} \sum_{m \in M} q_{rs}^{m+} \qquad (1)$$

subject to

$$G(\xi) \leq B \qquad (2)$$

$$o_r^+ = \sum_{s \in S} \sum_{m \in M} q_{rs}^{m+} \leq o_r^{\max} - o_r^0, \forall r \in R \qquad (3)$$

$$d_s^+ = \sum_{r \in R} \sum_{m \in M} q_{rs}^{m+} \leq d_s^{\max} - d_s^0, \forall s \in S \qquad (4)$$

$$\xi_n^m C_n^{m,\min} \leq \bar{c}_n^m \leq \xi_n^m C_n^{m,\max}, \forall m \in \bar{M}, n \in \bar{N}^m \qquad (5)$$

$$\sum_{r \in R} \sum_{s \in S} \sum_{p \in P_{rs}^m} \left( f_p^{m+} + f_p^{m0} \right) \Delta_{np} \leq \bar{c}_a^m, \forall m \in \bar{M}, n \in \bar{N}^m \qquad (6)$$

$$q_{rs}^+ \geq 0, \forall r \in R, s \in S \qquad (7)$$

$$\xi_n^m \in \{0,1\}, \forall m \in \bar{M}, n \in \bar{N}^m \qquad (8)$$

The objective of the upper level problem is to maximize the total number of trips generated. Equation (2) is the budget constraint, where the cost function $G(\xi)$ is generally assumed to be non-negative, increasing and differentiable (Yang and Bell, 1998). Equations (3) and (4), respectively, set the upper and lower bounds for the number of trips generated at each origin and destination node. Equation (5) is the capacity constraint for building transfer infrastructure. If a transfer location is determined, i.e., $\xi_n^m = 1$, then the capacity to build should be within $\left[ C_n^{m,\min}, C_n^{m,\max} \right]$. Equation (6) is the capacity constraint for the transfer node. The capacity constraint applies for the travel mode, in which the transfer flow is strictly bounded by the infrastructure capacity, i.e., number of parking slots. Equations (7) and (8) are definitional constraints for decision variables.

*Lower Level Problem*

The lower level model extends the combined trip distribution/assignment with variable costs model (Yang et al., 2000 for details) by introducing an additional modal choice component.

$$\min_{\mathbf{f}} Z = z_1 + z_2 + z_3 + z_4 + z_5 \qquad (9)$$

where

$$z_1 = \frac{1}{\theta} \sum_{m \in M} \sum_{p \in P_{rs}^m} f_p^m \left( \ln f_p^m - 1 \right) + \sum_{m \in M} \sum_{a \in A^m} \int_0^{v_a^m} t_a(x) dx \qquad (10)$$

$$z_2 = \sum_{m \in \bar{M}} \sum_{n \in \bar{N}} \int_0^{\bar{v}_n^m} \bar{t}_n^m(x) dx \qquad (11)$$

$$z_3 = \frac{1}{\gamma} \sum_{r \in R} \sum_{s \in S} \sum_{m \in M \cup \bar{M}} \left[ q_{rs}^{m+} \left( \ln \frac{q_{rs}^{m+}}{q_{rs}^+} - 1 \right) + q_{rs}^+ \right] \qquad (12)$$

$$z_4 = \frac{1}{\eta} \sum_{r \in R} \sum_{s \in S} q_{rs}^+ \left( \ln q_{rs}^+ - \ln o_r^+ \right) \qquad (13)$$

$$z_5 = -\sum_{s \in S} \int_0^{\sum_{r \in R} q_{rs}^+} h_s(x) dx \qquad (14)$$

Subject to:

$$\sum_{p \in P_{rs}^m} f_p^{m0} = q_{rs}^{m0}, \forall r \in R, s \in S, m \in M \cup \bar{M} \qquad (15)$$

$$\sum_{p \in P_{rs}^m} f_p^{m+} = q_{rs}^{m+}, \forall r \in R, s \in S, m \in M \cup \bar{M} \qquad (16)$$

$$q_{rs}^{+} = \sum_{m \in M \cup \bar{M}} q_{rs}^{m+}, \forall r \in R, s \in S \tag{17}$$

$$v_a^m = \sum_{r \in R} \sum_{s \in S} \sum_{m' \in M \cup \bar{M}} \sum_{p \in P_{rs}^{m'}} \left( f_p^{m'+} + f_p^{m'0} \right) \delta_{ap}, \forall a \in A, m \in M \tag{18}$$

$$\bar{v}_n^m = \sum_{r \in R} \sum_{s \in S} \sum_{p \in P_{rs}^m} \left( f_p^{m+} + f_p^{m0} \right), \forall m \in \bar{M}, n \in \bar{N}^m \tag{19}$$

$$f_p^{m+}, f_p^{m0} \geq 0, \forall m \in M \cup \bar{M}, p \in \bigcup_{r \in R, s \in S} P_{rs}^m \tag{20}$$

The objective function (9) of the lower level problem contains five elements defined via equations (10) - (14). $z_1$ and $z_2$, together, replicates the mathematical formulation for the stochastic user equilibrium problem, where $z_1$ and $z_2$ capture the cost of the links on each subnetwork and that at each transfer node. $z_3$ and $z_4$ respectively represent the entropy functions on mode choice and destination choice behavior. $z_5$ is the inverse of the elastic demand function. Equations (15) and (16) are flow conservation constraints. Equation (17) computes the newly generated flow at for each OD pair. Equations (18) and (19) are definitional constraints for link flow and flow at a transfer node, respectively. Equation (20) is nonnegativity constraints for flow variables. By examining the KKT conditions of the lower level problem, it can be proved in the appendix that the demand distribution, travelers' mode choice, and route choice follows multinomial logit distribution.

**SOLUTION ALGORITHM**

It is well known that a general network design problem is NP-hard and extremely difficult to solve. Therefore, various heuristic, metaheuristic or matheuristic methods have been proposed in the literature to solve the network design problem. In this study, we develop a matheurisc algorithm to solve the model. Due to the space limit, we are not able to provide the detail of the algorithm. In short, the algorithm integrates the Genetic Algorithm and a successive linear programming method developed in (Yang et al., 2000). In the GA, a chromosome represents the solutions associated with the transfer location and capacity, i.e., $\xi_n^m$ and $\bar{c}_n^m$ in the upper level problem. Given these variables, the successive linear programming method is applied to solve the resultant bilevel model and the objective value of the bilevel model is adopted as the fitness value of the chromosome. The algorithm terminates once a predefined maximum number of iterations is reached.

**NUMERICAL STUDY**

The numerical studies are designed to illustrate the properties of the model and demonstrate two Braess-like Paradox phenomena in a multimodal transport network. The first one states that under fixed demand, constructing more parking space to stimulate the usage of P+R service could deteriorate the system performance, measured by the total passengers' travel time. The second one states that under variable demand, increasing the parking capacity for the P+R to attract the usage may fail, represented by the decline in its modal share. The two phenomena, respectively, represent providing P+R services in a well-developed area, where the total travel demand stops growing, and in a developing area, where potential maximum travel demand has not been reached. The Paradox phenomena justify the necessity of optimising the capacity of parking space for the P+R services in terms of reducing the total travel time. Without further specified, the model is coded by Matlab 2018b and the lower level problem is solved by the fmincon function.

## Braess-like Paradox under Fixed Demand

*Occurrence of the Paradox*

The first experiment is conducted to demonstrate that introducing P+R services could induce a Braess-like paradox in a multimodal transport network under fixed OD demand. We construct a small network, which is similar to the classic Braess network, as shown in Figure 2, and consider 2000 travelers between nodes *O* and *D*. The network contains four nodes, where node *B* represents a metro station and node *A* denotes a potential parking area close to node B. There are four links and their travel modes and link performance functions are listed in the figure. Figure 2(a) and 2(b), respectively, represent the network before and after constructing parking slots at node *A* to provide the travel mode of P+R. In the before scenario, there are two travel modes and each mode has one only one path as shown in Figure 2(a), i.e., using private vehicle via path *O-A-D* and using metro via path *O-B-D*. In the after scenario, when a parking area is constructed at node *A*, the P+R mode is available to the travellers via path *O-A-B-D*.

**Link performance functions**

$$t_1^{car} = 4 + \left(\frac{v_1^{car}}{500}\right)^2, \quad t_2^{car} = 43 + \left(\frac{v_2^{car}}{1000}\right)^4, \quad t_3^{walk} = 25, \quad t_4^{metro} = 25 + \frac{v_4^{metro}}{500}$$

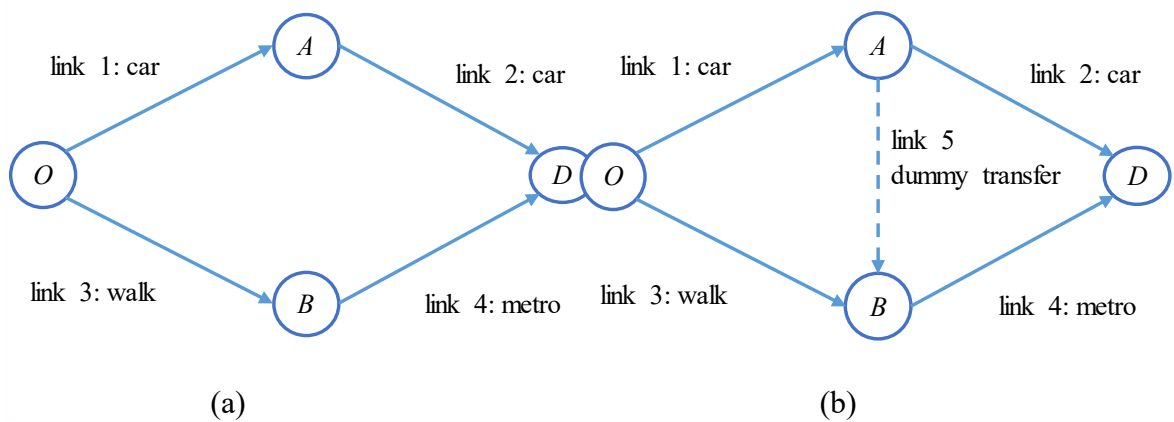

**Figure 2** Network for illustrating the Braess-like paradox phenomena

TABLE 1 summarizes total travel time and flow distribution obtained from the before and after scenarios. It shows that the total travel time in the after scenario is higher than that in the before scenario. This resembles the Braess paradox, which states that adding a new link deteriorates the network performance measured by the total travel time due to travelers' selfish route choice behavior. In our example, this is because the travel time of the P+R mode (path *O-A-B-D*) is much less than the other two modes, thus it attracts most of the travellers, leading the increase in the total travel time of the network.

**TABLE 1** Occurrence of the Braess-like paradox

| | Flow on Path 1 | Flow on Path 2 | Flow on Path 3 | Total travel time |
|---|---|---|---|---|

|         | O-A-D (Car) | O-B-D (Metro) | O-A-B-D (P+R) |        |
|---------|-------------|---------------|---------------|--------|
| Before  | 755         | 1245          | -             | 102790 |
| After   | 0           | 184           | 1816          | 102910 |

*Effect of $\theta$ on the Occurrence of the Paradox*

This section illustrates the effect of $\theta$ on the occurrence of the paradox. $\theta$ is a scaling parameter that measures travelers' perception on the randomness of the travel time. When its value increases to infinity, the result of the stochastic assignment model is close to that of the user equilibrium model. In this experiment, the value of $\theta$ is varied from 0.1 to 0.9. The resultant total travel times of the before and after scenarios are plotted in Figure 3(a), denoted by black and red curves, respectively. As shown in Figure 3(a), in the before scenario, the total travel time is stable. This is because the demand is distributed between the two paths in a way such that the travel times of the two paths are the same, despite the value of $\theta$. In contrast, in the after scenario, the total travel time is monotonically increasing with the increase of $\theta$ value. When the value of $\theta$ grows above 0.78, the total travel time in the after scenario is larger than that in the before scenario, manifesting the occurrence of the paradox. The result indicates that when the route choice behavior tends to follow user equilibrium (i.e., $\theta$ is large), it is more likely that the Braess-like paradox would happen.

*Effect of Capacity Constraints on the Occurrence of the Paradox*

This section examines the effect of capacity constraints of the parking space on the occurrence of the paradox, whiles it is ignored in previous examples by assuming that there is enough parking space. In the test, the capacity is varied from 100 to 2000 and the total travel times of the before and after scenarios at $\theta = 0.9$ are plotted in Figure 3(b). Interestingly, it is noticed that the paradox could occur either at a low capacity level (below 115) or a high capacity level (above 1824). The total travel time drops to bottom at when capacity is equal to 1300. The results justify the necessity of optimising the capacity of parking space for the P+R services in terms of reducing the total travel time.

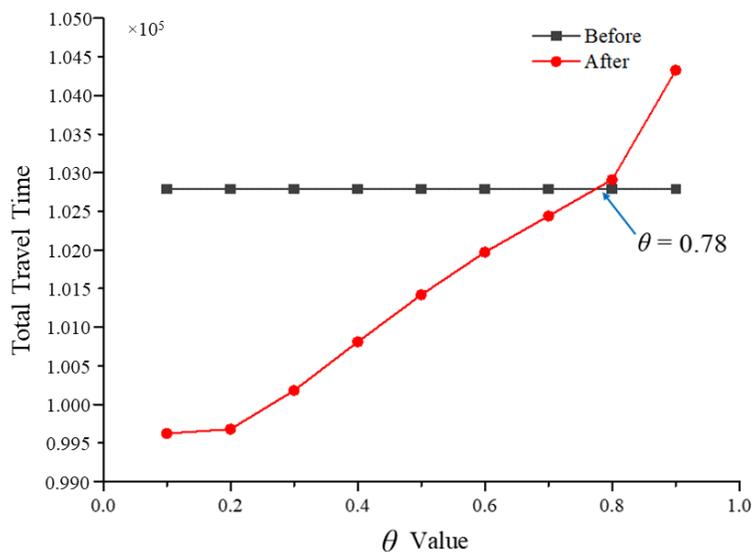

(a) Changes in the total travel time with the increase in $\theta$

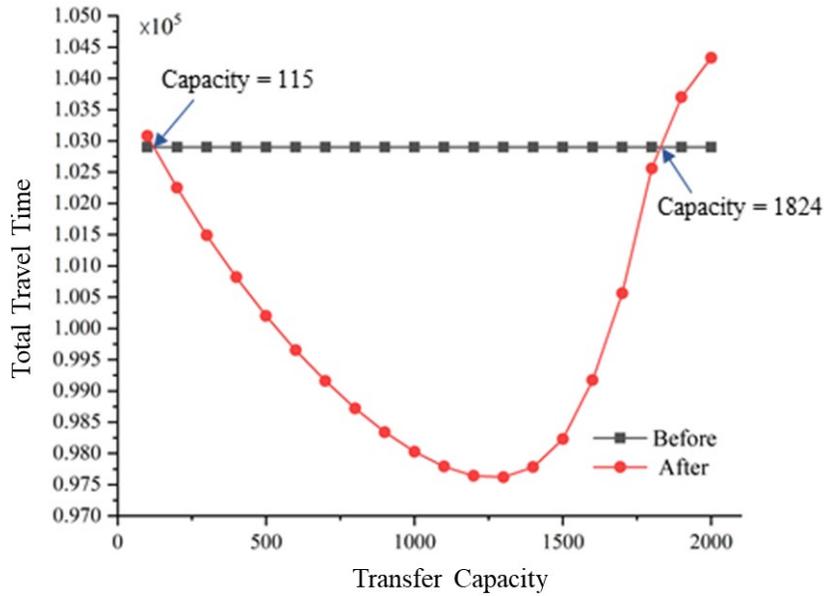

(b) Changes in the total travel time with the increase in the transfer capacity
**Figure 3** Effects $\theta$ and transfer capacity on the occurrence of the Braess-like paradox

**Paradox under Variable Demand**

In the experiments, we adopted the network in Figure 2 and varied the parking capacity from 100 to 2000. The existing and maximum demand are set as $o^0 = 100$ and $o^{max} = 2000$, respectively. The scaling parameters in the logit model is set at $\theta = 0.1$.

Under variable demand, the total travel time may not be a consistent measure due to the changes in the travel demand. Therefore, we examine the proportion of the travelers using the P+R with the increase in the parking capacity and plot the results in Figure 4. An interesting phenomenon is observed that expanding parking capacity could reduce the proportion of travelers using the P+R services. This indicates that adding more parking spaces does not make the P+R services more attractive to the newly generated travel demand. This is considered as a paradox, since it contracts with our intuition that increasing the P+R capacity would increase the attractiveness of such a service, which is in line with spirit of the Braess paradox and can be explained as a result of travelers selfish route choice behavior. Nevertheless, despite the reduction in the share of the P+R services, the market share of metro services still grows mildly. In terms of the path cost, as expected, it increases with in the increasing the in transfer capacity, since higher transfer capacity permits more travel demand to be generated in the network and induces higher congestion cost.

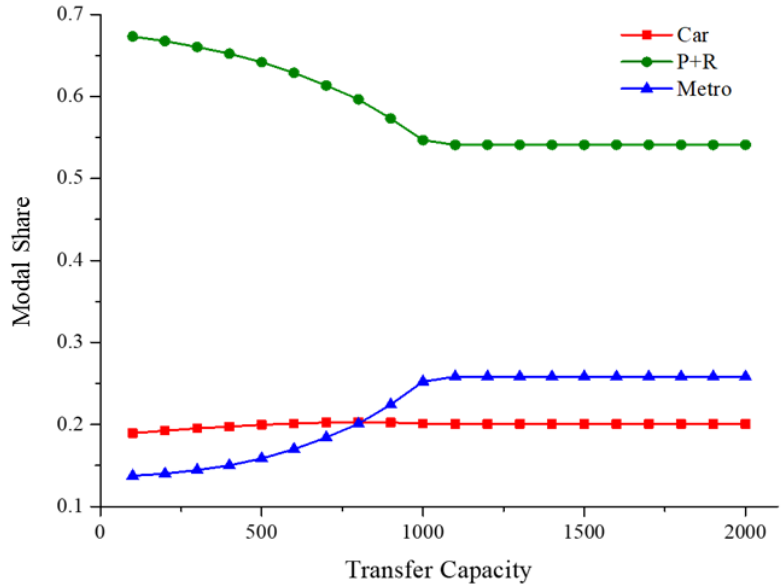

(a) Modal Share

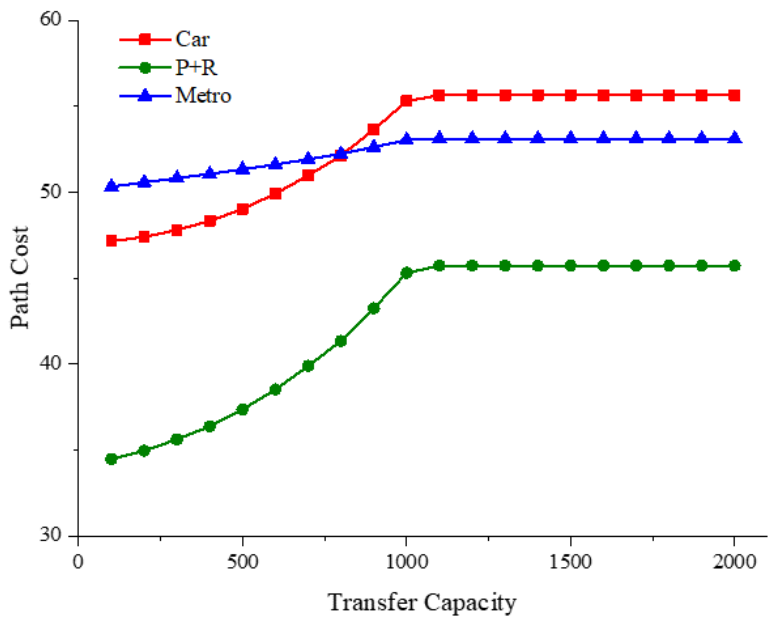

(b) Path Cost

**Figure 4** Paradox under variable demand

**Optimal Design of Multimodal Transfer Capacity**

This experiment is conducted to illustrate the design of optimal transfer capacity under the budget constraint in a multimodal transport network. We consider a more general network containing twelve links and seven nodes as shown in Figure 5. There are two transfer infrastructure options, one is to build parking spaces at node 7 to facilitate park-and-ride and the other is to build docking slots at node 5 to promote bike-and-ride. There are two OD pairs in this network, i.e. O-D pair (1,4) and O-D pair (6,4) with the existing demand of 300 and 500 respectively. We set the potential capacities of the B+R node and the P+R node from 300 to 1500 and 400 to 800 respectively. The budget constraint is 22,500,000 while the construction cost of each bike parking space and each car parking space is 12,500 and 25,000, respectively. The detailed route data and link performance functions are listed in TABLE 2.

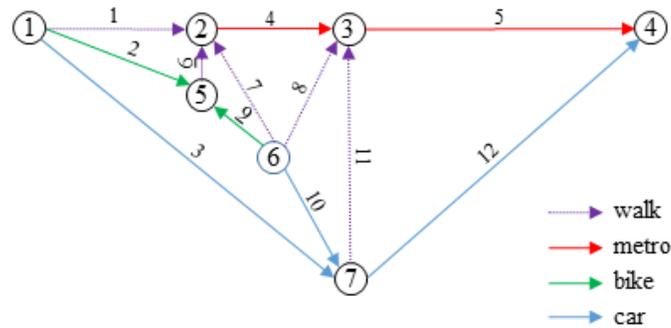

**Figure 5** The Transfer Nodes Design Network

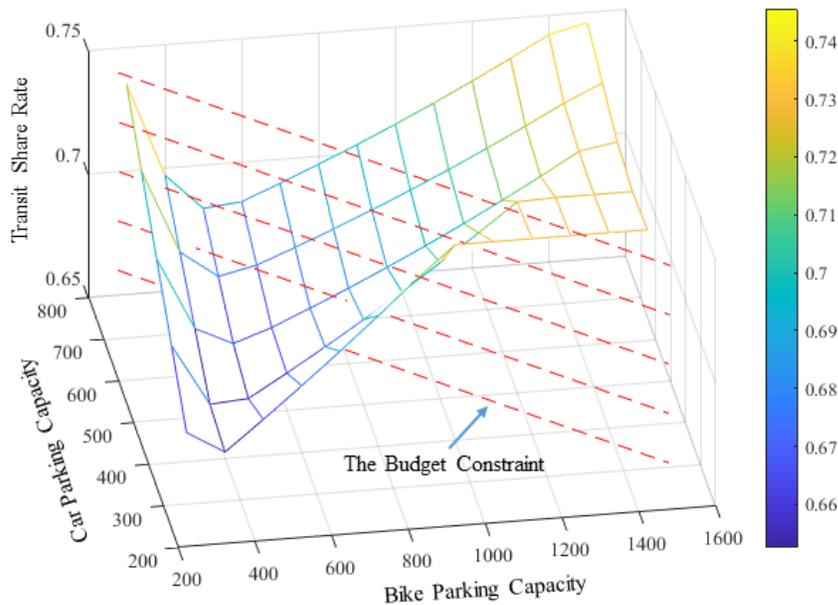

**Figure 6** Effect of transfer capacity on the share of sustainable transport mode under budget constraint

With the spread of sustainable design concept, we set the transit modes (metro, B+R and P+R) share rate as another reference value to see the optimal design of the transfer capacity. The share rates of transit modes resulted from the car parking capacity and bike parking capacity variation are shown in Figure 6. Both the enhancements of the bike parking capacity and the car parking capacity will firstly cause a decrease in the transit share rate and then an increase. Finally, the transit share rate will keep stable at 0.74 when the transfer capacity is large enough. As shown in Figure 6, with the budget constraint, the optimal transfer capacity is 400 bike parking spaces and 700 car parking spaces or 900 bike parking spaces and 450 car parking spaces. This result suggests it is better to distribute the budget to build to two transfer stations instead of establishing a large transfer center at either of the nodes.

**TABLE 2** Input data

| Link | Link Cost | Route | Nodes |
|---|---|---|---|
| Walk | | Car | |
| 1 | 20 | 2 | 1-7-4 |
| 6 | 2 | 6 | 6-7-4 |
| 7 | 12 | Metro | |
| 8 | 25 | 1 | 1-2-3-4 |
| 11 | 3 | 5 | 6-3-4 |
| Car | | 7 | 6-2-3-4 |
| 3 | $5+\left(\dfrac{v_3}{500}\right)^2$ | B+R | |
| 10 | $4+\left(\dfrac{v_{10}}{500}\right)^2$ | 3 | 1-5-2-3-4 |
| 12 | $40+\left(\dfrac{v_{12}}{1000}\right)^4$ | 8 | 6-5-2-3-4 |
| Metro | | P+R | |
| 4 | $3+\dfrac{v_4}{500}$ | 4 | 1-7-3-4 |
| 5 | $21+\dfrac{v_5}{500}$ | 9 | 6-7-3-4 |
| Bike | | Origins | Maximum Demand |
| 2 | 5 | $O_1$ | 2500 |
| 9 | 3 | $O_6$ | 2500 |

**CONCLUSION**

In view of the booming in the multimodal mobility, the paper develops a bilevel model for the multimodal network design problem. The upper level problem is to simultaneously determine the location and capacity of transfer infrastructures to be built. The lower level problem is the combined trip distribution and assignment model subjected to capacity constraints. It has been proved that at optimality the trip distribution and commuters' route choice behaviour depicted by the optimal solution of the lower level problem follow multinomial logit distribution. To solve the model, this paper employs a Genetic Algorithm to generate the solution associated with the transfer location and capacity and adopts a successive linear programming solution approach to solve the resultant bilevel model. Numerical studies were conducted to illustrate the optimal design of transfer capacity in a multimodal transport network and two Braess-like paradox phenomena, stating that constructing parking space to stimulate the usage of P+R service may induce higher total passengers' travel time and decline of the modal share of the P+R service.

This work opens various research directions. For example, 1) the lower level model assumes that the travelers' route choice following MNL model, which does not capture the nested structure of model choice and route choice as well as the overlap among paths. It will be more realistic to adopt other advanced behavior model; 2) the paradox phenomena are discussed in the context of constructing new infrastructure, it would be necessary to develop a methodology to detect whether existing infrastructure could trigger the paradox; 3) Existing lower level are coded by MATLAB and solved by its intrinsic function. It inhibits the application of solving large network. One undergoing research work is to develop a more efficient solution algorithm and an open source platform that favors large network applications.

## ACKNOWLEDGEMENT

This research is funded by the Scholarship Fund of China Scholarship Council (CSC No. 201806090216), the National Natural Science Foundation Council of China under Project No. 51638004, the Innovation Fund Denmark (IFD) under File No. 4109-00005 and the Jiangsu Graduate Research Innovation Program No. KYCX18_0135.

## AUTHOR CONTRIBUTIONS
The author confirm contribution to the paper as follows: Model Formulations (Jiao Ye and Yu Jiang); Coding(Jiao Ye and Yu Jiang); Paper writing (Yu Jiang and Jiao Ye); Paper Proofing (Jun Chen and Otto Anker Nielsen). All authors reviewed the results and approved the final version of the manuscript.

## APPENDIX A
This appendix shows that, at optimality, demand distribution, travelers' mode choice, and route choice obtained by solving the lower level problem follows multinomial logit distribution. This can be proved by examining the KKT conditions of the lower level problem. To begin with, we provide the Lagrange function of the lower level problem.

$$L = Z + \sum_{r \in R}\sum_{s \in S}\sum_{m \in M \cup \bar{M}} \alpha_{rs}^m (q_{rs}^{m0} - \sum_{p \in P_{rs}^m} f_p^{m0}) + \sum_{r \in R}\sum_{s \in S}\sum_{m \in M \cup \bar{M}} \beta_{rs}^m (q_{rs}^{m+} - \sum_{p \in P_{rs}^m} f_p^{m+})$$
$$+ \sum_{r \in R}\sum_{s \in S} \phi_{rs} (\sum_{m \in M \cup \bar{M}} q_{rs}^{m+} - q_{rs}^+) + \sum_{a \in A}\sum_{m \in M \cup \bar{M}} \lambda_a^m \left[ \sum_{r \in R}\sum_{s \in S}\sum_{m' \in M \cup \bar{M}}\sum_{p \in P_{rs}^{m'}} (f_p^{m'} + f_p^{m'0})\delta_{ap} - v_a^m \right] \quad (21)$$
$$+ \sum_{n \in N}\sum_{m \in M \cup \bar{M}} \mu_n^m \left[ \sum_{r \in R}\sum_{s \in S}\sum_{p \in P_{rs}^m} (f_p^{m+} + f_p^{m0}) - \bar{v}_n^m \right]$$

where $\alpha_{rs}^m$, $\beta_{rs}^m$, $\gamma_{rs}$, $\lambda_a^m$ and $\mu_n^m$, respectively, are the dual variables associated with Eqs. (15) - (19).

(1) To prove that the travelers' route choice follows multinomial logit distribution, we consider the KKT conditions w.r.t. to variable $f_p^m$, which is:

$$f_p^m \frac{\partial L}{\partial f_p^m} = 0 \quad \text{and} \quad \frac{\partial L}{\partial f_p^m} \geq 0, \ \forall p \in P_{rs}^m, \ m \in M \cup \bar{M} \quad (22)$$

$$\frac{\partial L}{\partial \alpha_{rs}^m} = 0, \ \frac{\partial L}{\partial \beta_{rs}^m} = 0, \ \frac{\partial L}{\partial \phi_{rs}} = 0, \ \frac{\partial L}{\partial \lambda_a^m} = 0, \ \frac{\partial L}{\partial \mu_n^m} = 0, \quad (23)$$
$$\forall a \in A, \ r \in R, \ s \in S, \ n \in N, \ m \in M \cup \bar{M}, \ p \in P_{rs}^m$$

where

$$\frac{\partial L}{\partial f_p^m} = \frac{\partial Z}{\partial f_p^m} + \frac{\partial}{\partial f_p^m} \sum_{n \in N} \sum_{m \in M \cup \bar{M}} \mu_n^m \left[ \sum_{r \in R} \sum_{s \in S} \sum_{p \in P_{rs}^m} (f_p^{m+} + f_p^{m0}) - \bar{v}_n^m \right]$$

$$= \frac{\partial}{\partial f_p^m} \frac{1}{\theta} \sum_{m \in M \cup \bar{M}} \sum_{p \in P_{rs}^m} f_p^m (\ln f_p^m - 1) + \frac{\partial}{\partial v_a^m} \sum_{m \in M \cup \bar{M}} \sum_{a \in A^m} \int_0^{v_a^m} t_a(x) dx \cdot \frac{\partial v_a^m}{\partial f_p^m} + \mu_n^m \quad (24)$$

$$= \frac{1}{\theta} \ln f_p^m + t_a(v_a^m) + \mu_n^m$$

Then, the KKT conditions can be transformed into

$$\ln f_p^m + t_a + \mu_n^m = 0 \Rightarrow \ln f_p^m = -\theta t_a - \theta \mu_n^m \Rightarrow \ln(f_p^m / f_{p'}^m) = -\theta(t_a - t_a') - \theta(\mu_n^m - \mu_n^{m'}) \quad (25)$$

which is consistent with the mathematical formulation of the stochastic user equilibrium.

(2) To prove that the travelers' mode choice follows multinomial logit distribution. We consider the KKT conditions w.r.t. to variable $q_{rs}^{m+}$, which is:

$$\frac{\partial L}{\partial q_{rs}^{m+}} = \frac{\partial Z}{\partial q_{rs}^{m+}} + \frac{\partial}{\partial q_{rs}^{m+}} \sum_{r \in R} \sum_{s \in S} \sum_{m \in M \cup \bar{M}} \beta_{rs}^m (q_{rs}^{m+} - \sum_{p \in P_{rs}^m} f_p^{m+})$$

$$+ \frac{\partial}{\partial q_{rs}^{m+}} \sum_{r \in R} \sum_{s \in S} \phi_{rs} (\sum_{m \in M \cup \bar{M}} q_{rs}^{m+} - q_{rs}^+)$$

$$= \frac{\partial}{\partial q_{rs}^{m+}} \frac{1}{\gamma} \sum_{r \in R} \sum_{s \in S} \sum_{m \in M} \left[ q_{rs}^{m+} (\ln \frac{q_{rs}^{m+}}{q_{rs}^+} - 1) + q_{rs}^+ \right] + \beta_{rs}^m + \phi_{rs} \quad (26)$$

$$= \frac{1}{\gamma} \left[ (\ln \frac{q_{rs}^{m+}}{q_{rs}^+} - 1) + q_{rs}^{m+} \frac{1}{q_{rs}^{m+}} \right] + \beta_{rs}^m + \phi_{rs}$$

$$= \frac{1}{\gamma} \ln \frac{q_{rs}^{m+}}{q_{rs}^+} + \beta_{rs}^m + \phi_{rs}$$

when it attains its minimum, the following condition holds

$$\frac{\partial L}{\partial q_{rs}^{m+}} = 0 \quad (27)$$

By substituting Equation (26) into (27), we have

$$q_{rs}^{m+} = q_{rs}^+ e^{-\gamma(\beta_{rs}^m + \phi_{rs})} \quad (28)$$

$$\sum_{m \in M \cup \bar{M}} q_{rs}^{m+} e^{-\gamma(\beta_{rs}^m + \phi_{rs})} = q_{rs}^+ \sum_{m \in M \cup \bar{M}} e^{-\gamma(\beta_{rs}^m + \phi_{rs})} = q_{rs}^+ \quad (29)$$

$$\sum_{m \in M \cup \bar{M}} e^{-\gamma(\beta_{rs}^m + \phi_{rs})} = 1 \quad (30)$$

$$q_{rs}^{m+} = q_{rs}^+ \frac{e^{-\gamma(\beta_{rs}^m + \phi_{rs})}}{\sum_{m \in M \cup \bar{M}} e^{-\gamma(\beta_{rs}^m + \phi_{rs})}} \quad (31)$$

which is consistent with the multinomial logit model.

(3) To prove that the demand distribution follows multinomial logit distribution. We consider the derivative of the Lagrange function w.r.t. to variable $q_{rs}^+$, which is

$$\frac{\partial L}{\partial q_{rs}^+} = \frac{\partial Z}{\partial q_{rs}^+} + \frac{\partial}{\partial q_{rs}^+} \sum_{r \in R} \sum_{s \in S} \phi_{rs} ( \sum_{m \in M \cup \bar{M}} q_{rs}^{m+} - q_{rs}^+ )$$

$$= \frac{\partial}{\partial q_{rs}^+} \frac{1}{\eta} \sum_{r \in R} \sum_{s \in S} q_{rs}^+ \left( \ln q_{rs}^+ - \ln q_r^+ \right) - \phi_{rs}$$

$$= \frac{1}{\eta} \left[ \left( \ln q_{rs}^+ - \ln q_r^+ \right) + q_{rs}^+ \frac{1}{q_{rs}^+} \right] - \phi_{rs}$$

$$= \frac{1}{\eta} (\ln \frac{q_{rs}^+}{q_r^+} + 1) - \phi_{rs} \tag{32}$$

When it attains its minimum, we have,

$$\frac{\partial L}{\partial q_{rs}^+} = 0 \tag{33}$$

By substituting Eq.(32) into (33), we obtain,

$$q_{rs}^+ = q_r^+ e^{\eta \phi_{rs} - 1} \tag{34}$$

$$\sum_{s \in S} q_r^+ e^{\eta \phi_{rs} - 1} = q_r^+ \sum_{s \in S} e^{\eta \phi_{rs} - 1} = q_r^+ \tag{35}$$

$$\sum_{s \in S} e^{\eta \phi_{rs} - 1} = 1 \tag{36}$$

$$q_{rs}^+ = q_r^+ \frac{e^{\eta \phi_{rs} - 1}}{\sum_{s \in S} e^{\eta \phi_{rs} - 1}} \tag{37}$$

which is consistent with the multinomial logit model.